\documentclass[10pt]{article}
\usepackage{amsmath,amstext, amsxtra, amsthm, amscd, amsgen, amsbsy, amsopn, amsfonts,
latexsym, amssymb, amscd, graphicx}

\newtheorem{theorem}{Theorem}
\newtheorem{corollary}[theorem]{Corollary}

\newtheorem{surgerylemma}[theorem]{Surgery lemma}

\def\proofend{\hbox to 1em{\hss}\hfill $\blacksquare $\bigskip }

\def\Z{{\mathbb Z}}
\def\R{{\mathbb R}}

\def\C{{\mathbb C}}

\def\H{{\mathbb H}}

\begin{document}

\title{A note on the $\hat A$-genus for $\pi_2$-finite manifolds with $S^1$-symmetry}
\author{Manuel Amann and Anand Dessai\footnote{Anand Dessai was partially supported by SNF Grant No. 200021-117701.}
}
\maketitle

\noindent
The purpose of this note is to answer the question whether the $\hat A$-genus vanishes on $S^1$-manifolds with finite second homotopy group. This question is connected to the work of Hayde\'e and Rafael Herrera \cite{HH} on $12$-dimensional positive quaternionic K\"ahler manifolds. To explain this we begin with a short incomplete survey of the classification problem for positive quaternionic K\"ahler manifolds (QK-manifolds) with special focus on the $12$-dimensional case. We refer to the survey article \cite{Sa} of Salamon for more information on QK-manifolds and references.

The only known examples of positive QK-manifolds are the symmetric examples studied by Wolf. LeBrun and Salamon showed that up to homothety there are only finitely many positive QK-manifolds in any fixed dimension and they conjectured that any positive QK-manifold is symmetric.

One knows that any positive QK-manifold $M$ is simply connected and that the second homotopy group $\pi _2(M)$ is trivial, isomorphic to $\Z $ or finite with $2$-torsion. In the first two cases $M$ is homothetic to the quaternionic projective space $\H P^n$ or the complex Grassmannian $Gr_2(\C ^{n+2})=U(n+2)/(U(n)\times U(2))$, respectively. There are symmetric examples, e.g. the Grassmannian $Gr_4(\R ^{n+4})=SO(n+4)/(SO(n)\times SO(4))$, which realize the third case. The question remains whether there exist non-symmetric positive QK-manifolds with finite second homotopy group.

The LeBrun-Salamon conjecture has been proved by Hitchin, Poon-Salamon and LeBrun-Salamon in dimension $\leq 8$. Hayde\'e and Rafael Herrera \cite{HH} showed that any $12$-dimensional positive QK-manifold $M$ is symmetric if the $\hat A$-genus of $M$ vanishes. If $M$ is a spin manifold this condition is always fulfilled by a classical result of Lichnerowicz since a positive QK-manifold has positive scalar curvature. One also knows that $\hat A(M)$ vanishes on the symmetric examples with finite second homotopy group (see \cite{BoHiII}, Th. 23.3).
Atiyah and Hirzebruch \cite{AtHi} showed that the $\hat A$-genus vanishes on spin manifolds with smooth effective $S^1$-action.

In \cite{HH} Hayde\'e and Rafael Herrera offered a proof for the vanishing of the $\hat A$-genus on any $\pi_2$-finite manifold with smooth effective $S^1$-action. Since one knows from the work of Salamon that the dimension of the isometry group of a $12$-dimensional positive QK-manifold is at least $5$ this would lead to a proof of the LeBrun-Salamon conjecture in this dimension.

The argument in \cite{HH} essentially consists of three parts. In the first part Hayde\'e and Rafael Herrera argue that any smooth $S^1$-action on a $\pi_2$-finite manifold is of even or odd type (this condition means that the sum of rotation numbers at the $S^1$-fixed points is always even or always odd). Then they argue that the proof of Bott-Taubes \cite{BoTa} for the rigidity of the elliptic genus may be adapted to non-spin manifolds if the $S^1$-action is of even or odd type. Finally they use an argument of Hirzebruch-Slodowy \cite{HiSl} to derive the vanishing of the $\hat A$-genus from the rigidity of the elliptic genus.

Unfortunately, the first part of their argument cannot be correct. In fact, as was noticed by the first named author, there are $S^1$-actions on the Grassmannian $Gr_4(\R ^{n+4})$ for any odd $n\geq 3$ which are neither even nor odd. For example, the $12$-dimensional Grassmannian $Gr_4(\R ^7)$ admits an $S^1$-action such that the fixed point components of the corresponding involution are of dimension $4$ and $6$ (the components are diffeomorphic to $S^4$ and $Gr_2(\R ^5)=SO(5)/(SO(3)\times SO(2))$ and both contain $S^1$-fixed points). However, for odd $n\geq 3$, $Gr_4(\R ^{n+4})$  is a non-spin positive QK-manifold with finite second homotopy group. The error in \cite{HH} can be traced back to an application of a result of Bredon on the representations at different fixed points which requires that $\pi _2(M)$ {\bf and}  $\pi _4(M)$ are finite (see the paragraph after Th. 4 in \cite{HH}).

This prompts the question whether one can prove the vanishing of the $\hat A$-genus on $\pi_2$-finite manifolds with smooth effective $S^1$-action by other means. The purpose of this note is to answer this question in the negative. More precisely, we will construct counterexamples in each dimension $4k\geq 8$ (in dimension $4$ the $\hat A$-genus does vanish on a simply connected $\pi_2$-finite manifold since it is a multiple of the signature). Our construction is a straightforward adaption of the classical elementary surgery theory (see \cite{Br}, Chapter IV) to the equivariant setting.

\begin{surgerylemma} Let $G$ be a compact Lie group and let $M$ be a smooth simply connected $G$-manifold. Suppose the fixed point manifold $M^G$ contains a submanifold $N$ of dimension $\geq 5$ such that the inclusion map $N\hookrightarrow M$ is $2$-connected. Then $M$ is $G$-equivariantly bordant to a simply connected $G$-manifold $M^\prime$ with $\pi _2(M^\prime)\subset \Z/2\Z$.\end{surgerylemma}

\noindent {\bf Proof:} Let $f:M\to BSO$ be a classifying map for the stable normal bundle of $M$. We fix a finite set of generators for the kernel of $f_*:\pi _2(M)\to \pi _2(BSO)\cong \Z /2\Z$. Since the inclusion map $N\hookrightarrow M$ is $2$-connected and $\dim N\geq 5$ we may represent these generators by disjointly embedded $2$-spheres in $N$. By construction the normal bundle in $M$ of each such $2$-sphere is trivial as a non-equivariant bundle and equivariantly diffeomorphic to a $G$-equivariant vector bundle over the trivial $G$-space $S^2$. For each embedded $2$-sphere we identify the normal bundle $G$-equivariantly with a tubular neighborhood of the sphere and perform $G$-equivariant surgery for all of these $2$-spheres. The result of the surgery is a simply connected $G$-manifold $M^\prime$ with $\pi _2(M^\prime)\subset \Z/2\Z$ (if $M$ is a spin manifold then $M^\prime$ is actually $2$-connected).\proofend

\begin{corollary}
 For any $k>1$ there exists a smooth simply connected $4k$-dimen\-sio\-nal $\pi_2$-finite manifold $M_{4k}$ with smooth effective $S^1$-action and $\hat A(M_{4k})\neq 0$.
\end{corollary}

\noindent {\bf Proof:} We begin with some linear effective $S^1$-action on the complex projective space $\C P^{2k}$ such that the fixed point manifold $M^{S^1}$ contains a component $N$ diffeomorphic to $\C P^l$ for some $l\geq 3$. Since $N\hookrightarrow \C P^{2k}$ is $2$-connected, the manifold $\C P^{2k}$ is $S^1$-equivariantly bordant to a simply connected $S^1$-manifold $M^\prime$ with $\pi _2(M^\prime)$ finite by the surgery lemma (in fact, $\pi _2(M^\prime)\cong \Z /2\Z $ since $\C P^{2k}$ is not a spin manifold). It is well-known that the $\hat A$-genus does not vanish on $\C P^{2k}$. Since $M^\prime $ is bordant to $\C P^{2k}$ we get $\hat A(M^\prime )=\hat A(\C P^{2k})\neq 0$.\proofend

It is straightforward to produce examples with much larger symmetry using the construction above. We leave the details to the reader.

It remains a challenging task to determine whether the $\hat A$-genus vanishes on $\pi _2$-finite positive QK-manifolds as predicted by the LeBrun-Salamon conjecture.

\end{document}